\documentclass[12pt]{amsart}
\usepackage{amsmath}
\usepackage{amsfonts}
\usepackage[T1]{fontenc}

\newtheorem{theorem}{Theorem}[section]
\newtheorem{lemma}[theorem]{Lemma}
\newtheorem{proposition}[theorem]{Proposition}
\newtheorem{corollary}[theorem]{Corollary}

\theoremstyle{definition}
\newtheorem{definition}[theorem]{Definition}
\newtheorem{notation}[theorem]{Notation}

\theoremstyle{remark}
\newtheorem{remark}[theorem]{Remark}

\newcommand{\bl}{\begin{lemma}}
\newcommand{\el}{\end{lemma}}
\newcommand{\bpr}{\begin{proposition}}
\newcommand{\epr}{\end{proposition}}
\newcommand{\bd}{\begin{definition}}
\newcommand{\ed}{\end{definition}}
\newcommand{\br}{\begin{remark}}
\newcommand{\er}{\end{remark}}
\newcommand{\bt}{\begin{theorem}}
\newcommand{\et}{\end{theorem}}
\newcommand{\bc}{\begin{corollary}}
\newcommand{\ec}{\end{corollary}}
\newcommand{\bn}{\begin{notation}}
\newcommand{\en}{\end{notation}}
\newcommand{\bp}{\begin{proof}}
\newcommand{\ep}{\end{proof}}
\newcommand{\be}{\begin{equation}}
\newcommand{\ee}{\end{equation}}
\newcommand{\nd}{\noindent}
\newcommand{\N}{{\mathbb{N}}}
\newcommand{\cn}{{\mathcal{N}}}
\newcommand{\cm}{\mathcal{M}}
\newcommand{\supp}{{\overline{supp}}}
\newcommand{\abs}[1]{\lvert#1\rvert}
\newcommand{\abd}[1]{\lvert#1\rvert^*}
\newcommand{\nrm}[1]{\|#1\|}
\newcommand{\nrd}[1]{\|#1\|^*}

\newcommand{\al}{\alpha}
\newcommand{\vep}{\varepsilon}

\newcommand{\bnum}{\begin{enumerate}}
\newcommand{\enum}{\end{enumerate}}
\newcommand{\ct}{{\mathcal{T}}}
\newcommand{\cs}{{\mathcal{S}}}
\newcommand{\cf}{{\mathcal{F}}}
\newcommand{\fn}{{[\N]^{<\infty}}}

\begin{document}

\begin{center}
\large{\textbf{Note on distortion and Bourgain $\ell_1$-index}}

\

\textsc{Anna Maria Pelczar}

\

\end{center}



{\small \textsc{Abstract}. The relation between different notions
measuring proximity to $\ell_1$ and distortability of a Banach
space is studied. The main result states that a Banach space,
whose all subspaces have Bourgain $\ell_1$-index greater than
$\omega^\alpha$, $\al<\omega_1$, contains either an arbitrary
distortable subspace or an $\ell_1^\al$-asymptotic subspace.}

\section{Preliminaries}

The study of asymptotic properties and in particular complexity of
the family of copies of $\ell_1^n$ in Banach spaces is closely
related to investigating their distortability, cf.
\cite{mt,mr,otw}. The investigation of arbitrary distortion of
Banach spaces is concentrated mainly on $\ell_1$-asymptotic
spaces. The first tool measuring the way $\ell_1$ is represented
in a Banach space is provided by Bourgain $\ell_1$-index. Another
approach is given by higher order spreading models, studied
extensively in mixed and modified mixed Tsirelson spaces. The
$\ell_1$-asymptoticity of higher order of a Banach space can be
measured by constants introduced in \cite{otw}.

We present here an observation, in the spirit of the theorem of \cite{mt} recalled below,
relating distortability of a Banach space to the "proximity" to $\ell_1$ measured by the
tools presented above.

\bt \cite{mt}\label{mt} Let $X$ be a Banach space. Then $X$
contains either an arbitrarily distortable subspace or an $\ell_p$
($1\leq p<\infty$) or $c_0$-asymptotic subspace.\et

Our main result states that a Banach space with a basis, whose all
block subspaces have Bourgain $\ell_1$-block index greater than
$\omega^\al$, contains either an arbitrary distortable subspace or
a $\ell_1^\al$-asymptotic subspace. In particular a space
saturated with $\ell_1^\al$-spreading models generated by block
sequences contains either an arbitrary distortable subspace or an
$\ell_1^\al$-asymptotic subspace. We obtain as a corollary the
theorem presented in \cite{m}. Analogous results hold also in
$c_0$ case. As a corollary we obtain that the "stabilized" (with
respect to block subspaces) Bourgain $\ell_1$-block index of a
space with bounded distortion not containing $\ell_1$ is of the
form $\omega^{\omega^\gamma}$ for some non-limit
$\gamma<\omega_1$.

\

We recall now the basic definitions and standard notation. Let $X$
be a Banach space with a basis $(e_i)$. A \textit{support} of a
vector $x=\sum_i x_i e_i$ is the set $supp\ x =\{ i\in\N :\
x_i\not =0\}$. We shall use also \textit{interval support} of a
vector $x\in X$ - the smallest interval in $\N$ containing support
of $x$ and denote it by $\supp\ x$.

Given any $x=\sum_i x_ie_i$ and finite $E\subset\N$ put
$Ex=\sum_{i\in E}x_ia_i$. We write $x<y$ for vectors $x,y\in X$,
if $\max (supp\ x)<\min (supp\ y)$. A \textit{block sequence} is
any sequence $(x_i)\subset X$ satisfying $x_{1}<x_{2}<\dots$, a
\textit{block subspace} of $X$ - any closed subspace spanned by an
infinite block sequence. If $Y$ is a block subspace of $X$ spanned
by a block sequence $(y_i)$ then $Y_n$, $n\in\N$, denotes the
"tail" subspace spanned by $(y_i)_{i\geq n}$ and $EY$, $E\subset
\N$, denotes the subspace spanned by $(y_i)_{i\in E}$.

A basic sequence $(x_1,\dots,x_k)$ in a Banach space is
$K$-\textit{equivalent} to the unit vector basis of
$k$-dimensional space $\ell_1$ (resp. $c_0$), for some $K\geq 1$,
if for any scalars $a_1,\dots,a_k$ we have
$K\nrm{a_1x_1+\dots+a_kx_k}\geq \abs{a_1}+\dots+\abs{a_k}$ (resp.
$\nrm{a_1x_1+\dots+a_kx_k}\leq
K\max\{\abs{a_1},\dots,\abs{a_k}\}$).

\bd A Banach space $(X, \nrm{\cdot})$ is
$\lambda-$\textit{distortable}, for $\lambda>1$, if there is an
equivalent norm $|\cdot|$ on $X$ such that for any infinite
dimensional subspace $Y$ of $X$
 $$\sup\left\{\frac{|x|}{|y|}:\ x,y\in Y,\ \nrm{x}=\nrm{y}=1
\right\}\geq\lambda$$ A Banach space $X$ is \textit{arbitrary
distortable}, if it is $\lambda-$distortable for any $\lambda>1$.

A Banach space $X$ has $D$-\textit{bounded distortion}, if for any
equivalent norm $\abs{\ \cdot\ }$ and any infinite dimensional
subspace $Y$ of $X$ there is a further infinite dimensional
subspace $Z$ of $Y$ such that $\abs{x}/\nrm{x}\leq
D\abs{y}/\nrm{y}$ for any non-zero $x,y\in Z$. A Banach space has
\textit{bounded distortion} if it has $D$-bounded distortion for
some $D\geq 1$.\ed

Notice that any Banach space $X$ contains either an arbitrary
distortable subspace or a subspace with bounded distortion.

\

Given any $M\subset\N$ by $[M]^{<\infty}$ denote the family of
finite subsets of $M$. A family $\cf\subset \fn$ is
\textit{regular}, if it is \textit{hereditary}, i.e. for any
$G\subset F$, $F\in \cf$ also $G\in \cf$, \textit{spreading}, i.e.
for any integers $n_1<\dots<n_k$ and $m_1<\dots<m_k$ with $n_i\leq
m_i$, $i=1,\dots, k$, if $(n_1,\dots,n_k)\in \cf$ then also
$(m_1,\dots,m_k)\in \cf$, and \textit{compact} in the product
topology of $2^\N$. If $\cf\subset\fn$ is compact, let $\cf'$
denote the set of limit points of $\cf$. Define inductively
families $\cf^{(\al)}$, $\al$ ordinal, by putting $\cf^{(0)}=\cf$,
$\cf^{(\al+1)}=(\cf^{(\al)})'$ for any $\al$ ordinal and
$\cf^{\al}=\bigcap_{\xi<\al}\cf^{(\xi)}$ for any $\al$ limit
ordinal. Define index $\iota(\cf)=\inf\{\al:\
\cf^{(\al+1)}=\emptyset\}$.

\

A \textit{tree on a set} $S$ is a subset of $\bigcup_{n=1}^\infty S^n$ such that
$(x_1,\dots,x_k)\in\ct$ whenever $(x_1,\dots,x_k,x_{k+1})\in\ct$, $k\in\N$. A tree $\ct$
is \textit{well-founded}, if there is no infinite sequence $(x_i)\subset S$ with
$(x_1,\dots,x_k)\in\ct$ for any $k\in\N$. Given a tree $\ct$ on $S$ put
$$D(\ct)=\{(x_1,\dots,x_k):\ (x_1,\dots,x_k,x)\in\ct \ \ \mathrm{for}\ \ \mathrm{some}\ \ x\in S\}$$
Inductively define trees $D^\al(\ct)$: $D^0(\ct)=\ct$,
$D^{\al+1}=D(D^\al(\ct))$ for $\al$ ordinal and
$D^\al(\ct)=\bigcap_{\xi<\al}D^\xi(\ct)$ for $\al$ limit ordinal.
The \textit{order} of a well-founded tree $\ct$ is given by
$o(\ct)=\inf\{\al:\ D^\al(\ct)=\emptyset \}$.

\

Let $X$ be a Banach space with a basis. A tree $\ct$ on $X$ is an
$\ell_1$-$K$-\textit{block tree on} $X$, $K\geq 1$, if any
$(x_1,\dots,x_k)\in \ct$ is a normalized block sequence
$K$-equivalent to the unit vector basis of $k$-dimensional space
$\ell_1$. An $\ell_1$-block tree on $X$ is an $\ell_1$-$K$-block
tree on $X$ for some $K\geq 1$.

Let $I_b(X,K)=\sup\{o(\ct):\ \ct$ is a $\ell_1$-$K$-block tree on
$X\}$, $K\geq 1$. The \textit{Bourgain $\ell_1$-block index of
$X$} is defined by $I_b(X)=\sup\{I_b(X,K):\ K\geq 1\}$.

\bt \cite{jo} \label{jo} Let $X$ be a Banach space with a basis not containing $\ell_1$.
Then $I_b(X)=\omega^\al$ for some $\al<\omega_1$ and $I_b(X)>I_b(X,K)$ for any $K\geq
1$.\et

\br Recall the close relation (\cite{jo}) between  $I_b(X)$ and
$I(X)$ - the original Bourgain $\ell_1$-index defined as block
index but by trees of not necessarily block sequences: for
$I(X)\geq \omega^\omega$ we have $I_b(X)=I(X)$, if
$I(X)=\omega^{n+1}$, then $I_b(X)=\omega^{n+1}$ or $\omega^n$,
$n\in\N$. \er

The \textit{generalized Schreier families} $(\cs_\al)_{\al<\omega_1}$ of finite subsets
of $\N$, introduced in \cite{aa}, are defined by the transfinite induction:
$$\cs_0=\{\{ n\}:\ n\in\N\}\cup\{\emptyset\}$$
Suppose the families $\cs_\xi$ are defined for all $\xi<\al$. If $\al=\beta+1$, put
$$\cs_\al=\left\{ F_1\cup\dots\cup F_m:\ m\in\N,\ F_1,\dots,F_m\in \cs_\beta,\ m\leq F_1<\dots<F_m\right\}$$
If $\al$ is a limit ordinal, choose $\al_n\nearrow \al$ and set
$$\cs_\al=\{F:\ F\in \cs_{\al_n}\ \mathrm{and}\ n\leq F\ \mathrm{for\ some}\ n\in\N\}$$

\nd It is well known that any family $\cs_\al$, $\al<\omega_1$, is regular with
$\iota(\cs_\al)=\omega^\al$, considered as a tree on $\N$ satisfies
$o(\cs_\al)=\omega^\al$, cf. \cite{aa}.

\

Fix $\al<\omega_1$. A finite sequence $(E_i)$ of subsets of $\N$ is
$\al-$\textit{admissible} (resp. $\al-$\textit{allowable}), if $E_1<E_2<\dots$ (resp.
$(E_i)$ pairwise disjoint) and $(\min E_i)\in \cs_\al$.

Let $X$ be a Banach space with a basis $(e_n)$, fix
$\al<\omega_1$. A finite sequence $(x_i)\subset X$ is
$\al$-admissible (resp. $\al$-allowable) with respect to the basis
$(e_n)$, if $( supp\ x_i)$ is $\al$-admissible (resp.
$\al$-allowable).

\bd\label{smas} Fix $\al<\omega_1$. Let $X$ be a Banach space with a basis $(e_n)$.

A normalized block sequence $(x_i)\subset X$ \textit{generates an
$\ell_1^\al$-spreading model} with constant $C\geq 1$, if for any
$F\in \cs_\al$ the sequence $(x_i)_{i\in F}$ is $C$-equivalent to
the unit vector basis of $\sharp F$-dimensional space $\ell_1$.

The space $X$ is $\ell_1^\al$-\textit{asymptotic} (resp.
$\ell_1^\al$-\textit{strongly asymptotic}) with constant $C\geq
1$, if for any $\al$-admissible (resp. $\al$-allowable) w.r.t.
$(e_n)$ sequence $(x_i)_{i=1}^k$ is $C$-equivalent to the unit
vector basis od $k$-dimensional space $\ell_1$.

\ed

Obviously $X$ is $\ell_1^\al$-asymptotic with constant $C$ iff all
normalized block sequences in $X$ generate $\ell_1^\al$-spreading
model with constant $C$. By the properties of $\cs_\al$'s, any
block subspace of an $\ell_1^\al$-asymptotic (resp. strongly
asymptotic) space with constant $C$ is also
$\ell_1^\al$-asymptotic (resp. strongly asymptotic) with the same
constant. The relations between Bourgain $\ell_1$-block index and
notions introduced above are described by

\bpr\label{relacja} Let $X$ be a Banach space with a basis. Fix $\al<\omega_1$.

If $X$ admits an $\ell_1^\al$-spreading model, then $I_b(X)>
\omega^\al$.

If $X$ is an $\ell_1^\al$-asymptotic space, then
$I_b(X)\geq\omega^{\al\omega}$.\epr

\bp The first part follows from Theorem \ref{jo} and the fact that
$o(\cs_\al)=\omega^\al$. The second part follows from the proof of Theorem 5.19,
\cite{jo}. We recall it briefly. For any $\cm,\cn\subset \fn$ put
$$\begin{array}{rl}\cm[\cn] =& \{  F_1\cup\dots\cup F_k:\ F_1,\dots,F_k\in
\cn,\ \ m_1\leq F_1<\dots<m_k\leq F_k \\ & \ \mathrm{for}\ \
\mathrm{some}\ \ (m_1,\dots,m_k)\in\cm,\ k\in\N\}\end{array}$$
 Put $[\cs_\al]^n=\cs_\al[\dots[\cs_\al]]$ ($n$ times). If $X$ is
$\ell_1^\al$-asymptotic with constant $C$, then any normalized
block sequence $(x_1,\dots,x_k)$ with $(\min (supp\ x_i))\in
[\cs_\al]^n$ is $C^n$-equivalent to the unit vector basis of
$k$-dimensional space $\ell_1$. Since $o([\cs_\al]^n)=\omega^{\al
n}$, \cite{aa}, therefore $I_b(X)>\omega^{\al n}$ for any
$n\in\N$, which ends the proof. \ep

\br The Definition \ref{smas} extends the well-known notions of
$\ell_1$-asym\-pto\-tic space (introduced in \cite{mt}) and
spreading model generated by a basic sequence. The higher order
$\ell_1$-spreading models were introduced in \cite{kn} and
investigated in particular in \cite{adm,lt2,m}. The constants
describing $\ell_1$-asymptoticity of higher order were introduced
and studied in \cite{otw}. The term $\ell_1$-asymptoticity of
higher order was explicitly introduced in \cite{g3}, where also a
criterium for arbitrary distortion in terms of $\ell_1$-spreading
models was given. The $\ell_p$-strongly asymptotic spaces were
introduced and studied in \cite{dfko}. Bourgain $\ell_1$-index and
$\ell_1$-block index of various spaces in relation to existence of
higher order spreading models, distortability and quasiminimality
were investigated in particular in \cite{jo,klmt,lt1,lt2}.\er

We shall need additional norms given by the $\ell_1$-asymptoticity
of the space:

\bd\label{norm} Let $U$ be a Banach space with a basis. Fix
$\al<\omega_1$. If $U$ is $\ell_1^\al$-asymptotic with constant
$C$, define an associated norm $\abs{\ \cdot\ }_\al$ on $U$ by
$$\abs{x}_\al=\sup\left\{\sum_{i=1}^k\nrm{E_ix}:\ E_1<\dots<E_k\ \al-\mathrm{admissible},\ k\in\N\right\},\ x\in
U$$
 If $U$ is $\ell_1^\al$-strongly asymptotic, in the definition of associated norm we use allowable
sequences instead of admissible ones. Clearly
$\nrm{\cdot}\leq\abs{\ \cdot\ }_\al\leq C\nrm{\cdot}$. \ed

Simpler versions of these norms were used to show arbitrary
distortion of the famous Schlumprecht space, the first Banach
space known to be arbitrary distortable, these norms also distort
some mixed and modified mixed Tsirelson spaces \cite{ad,adm,lt2}.

\section{Main result}

Now we present the main result, which shows that we can reverse the implication in Prop.
\ref{relacja} in spaces with bounded distortion.

\bt\label{main} Let $X$ be a Banach space with a basis. Fix
$\al<\omega_1$. Assume that $I_b(Y)>\omega^\al$ for any block
subspace $Y$ of $X$. Then $X$ contains either an arbitrary
distortable subspace or an $\ell_1^\al$-asymptotic subspace.

If, additionally, $X$ is $\ell_1^1$-strongly asymptotic, then $X$
contains either an arbitrary distortable subspace or an
$\ell_1^\al$-strongly asymptotic subspace.\et

By Prop. \ref{relacja} a Banach space admitting in any block
subspace $\ell_1^\al$-spreading models generated by normalized
block sequences satisfies the assumptions of the Theorem
\ref{main}. We have also the following Corollary, implying Theorem
2.1, \cite{m}:

\bc \label{col} Fix $1<\al<\omega_1$. Let a Banach space $X$ with
a basis admit for any $\xi<\al$ in every block subspace
$\ell_1^\xi$-spreading model generated by a normalized block
sequence with a universal constant $C\geq 1$. Then $X$ contains
either an arbitrary distortable subspace or an
$\ell_1^\al$-asymptotic subspace.

If, additionally, $X$ is $\ell_1^1$-strongly asymptotic, then $X$
contains either an arbitrary distortable subspace or an
$\ell_1^\al$-strongly asymptotic subspace.\ec

\bp Assume $X$ has no arbitrary distortable subspaces. If
$\al=\beta+1$ for some $\beta<\omega_1$, then by Theorem
\ref{main}, there is a $\ell_1^\beta$-asymptotic subspace $W$ with
some constant $C\geq 1$. By Prop. 3.2 \cite{otw} there is
$n_0\in\N$ such that $F\in \cs_\beta$ for any $n_0\leq F\in
\cs_1$. Thus $W_{n_0}$ is also $\ell_1^1$-asymptotic with constant
$C$, and therefore also $\ell_1^\al$-asymptotic (with constant
$C^2$).

If $\al$ is a limit ordinal, then by assumption
$I_b(Y)>I_b(Y,C)\geq\omega^\al$ for any block subspace $Y$ of $X$,
and Theorem \ref{main} ends the proof. The case of $\ell_1$-strong
asymptoticity follows analogously. \ep

\br By Lemma 6.5, \cite{jo}, (Remark 6.6 (iii)) the universal
constant $C$ (arbitrarily close to 1) in the assumption of the
Corollary \ref{col} is automatic for $\al=\omega^\gamma$, with
$\gamma$ limit ordinal.\er

\bc Let $X$ be a Banach space of bounded distortion with a basis,
not containing $\ell_1$. If $I_b(Y)=I_b(X)$ for any block subspace
$Y$ of $X$, then $I_b(X)=\omega^{\omega^\gamma}$ for some
non-limit $\gamma<\omega_1$.\ec

\bp Let $I_b(X)=\omega^\al$. For any $\beta<\al$, by Theorem
\ref{main}, $X$ has a $\ell_1^\beta$-asymptotic subspace, thus by
Prop. \ref{relacja}, $I_b(X)>\omega^{\beta 2}$. Hence for any
$\beta<\al$ also $\beta 2<\al$, thus $\al=\omega^\gamma$ for some
$\gamma<\omega_1$. By Remark 5.15 (iii), \cite{jo}, $\gamma$ is
not a limit ordinal. \ep

\br Observe that any Banach space $X$ has a block subspace $Y$
with $I_b(Z)=I_b(Y)$, for any block subspace $Z$ of $Y$. Indeed,
either $X$ contains $\ell_1$, or $I_b(X)<\omega_1$ (\cite{b}) and
we can use standard diagonalization.\er

\br\label{ex} We collect some known examples:

\bnum \item[(i)] $I_b(X)>\omega$ iff $1$ belongs to Krivine set of
$X$, i.e. $\ell_1$ is finitely (almost isometrically) represented
on block sequences in $X$.

\item[(ii)] For any $\al<\omega_1$ by Theorem 5.19, \cite{jo}, and
Prop. \ref{relacja}, any block subspace $Y$ of the Tsirelson type
space $T(\cs_\al, 1/2)$ (which is clearly $\ell_1^\al$-asymptotic)
satisfies $I_b(Y)=\omega^{\al\omega}$, thus $I_b(X)$ is of the
form $\omega^{\omega^{\gamma+1}}$.

\item[(iii)] By Theorem 4.2, \cite{ao}, mixed Tsirelson space
$X=T[(\cs_n,\theta_n)_n]$, $\theta_n\searrow 0$ contains no
$\ell_1^\omega$-asymptotic subspace. On the other hand it was
shown in \cite{lt2} that $I_b(Y)>\omega^\omega$ for any block
subspace $Y$ of $X$ iff any block subspace $Y$ admits
$\ell_1^\omega$-spreading model. In such a case $X$ is arbitrary
distortable. It holds in particular if $\lim\sqrt[n]{\theta_n}=1$.

\item[(iv)] In \cite{lt1} Bourgain $\ell_1$-block index of mixed
Tsirelson spaces was computed, as a consequence it is proved that
for any $\al$ not of the form $\omega^\gamma$, $\gamma$ limit
ordinal, there is a Banach space $X_\al$ with
$I_b(X_\al)=\omega^\al$. In particular it was proved that (with a
special choice of sequences in definition of Schreier families)
$I_b(T(\cs_{\beta_n},\theta_n)_n)$ is either $\omega^{\omega^\xi
2}$ or $\omega^{\omega^\xi}$, where $\beta_n\nearrow\omega^\xi$,
$\xi<\omega_1$ successor.  \enum\er

\nd \textit{Proof of Theorem \ref{main}}. We can assume that $X$
has a bimonotone basis. Assume $X$ contains no arbitrary
distortable subspaces, and pick a block subspace $Y$ of $X$ with
$D$-bounded distortion, for some $D\geq 1$. We restrict our
consideration to $Y$ and use the transfinite induction.

\

Idea of the proof of the first inductive step and limit case of
the second inductive step (the successor case is trivial) could be
described as follows: we consider equivalent norms, whose uniform
equivalence to the original norm would give asymptoticity of
desired order. We "glue" the norms on some special vectors
provided by high $\ell_1$-index of the space (Lemmas \ref{lem1}
and \ref{lem2}), using methods standard now in the study of
Tsirelson type spaces, and by the bounded distortion of a space we
obtain a uniform equivalence of these norms to the original one on
some subspace.

\

\nd\textsc{First inductive step}

The result for $\al=1$ follows from Theorem \ref{mt}, but we present here a shorter
proof, whose idea was used in the proof of Theorem \ref{mt} given in \cite{mr}, and whose
scheme will serve also in the next step. Define new equivalent norms on $Y$ as follows:
$$\nrm{y}_n=\sup\left\{\sum_{j=1}^n\nrm{E_jy}:\ E_1<\dots<E_n \ \mathrm{intervals}\right\},\ \ y\in Y,\ \ n\in\N$$
We recall a standard observation providing vectors "gluing" the
original norm and the new norms:

\bl\label{lem1} Let $U$ be a Banach space with a bimonotone basis.
Fix $n\in\N$ and assume $I_b(U)>\omega$. Then there is a vector
$x\in U$ with $1/2\leq\nrm{x}\leq \nrm{x}_n\leq 2$. \el

\nd \textit{Proof of Lemma \ref{lem1}}. By Remark \ref{ex} (i)
take a normalized block sequence $(x_i)_{i=1}^{n^2}\subset U$
which is 2-equivalent to the unit vector basis of
$n^2-$dimensional $\ell_1$ and put
$x=\frac{1}{n^2}\sum_{i=1}^{n^2}x_i$. Obviously $\nrm{x}\geq 1/2$.
Take any $E_1<\dots<E_n$ and put
$$I=\left\{ i:\ \supp\ x_i\
\mathrm{intersects}\ E_j\ \mathrm{and}\ E_{j+1}\ \mathrm{for \
some}\ j \right\}$$

\nd Since $\sharp I\leq n$ we have
$$\sum_{j=1}^n\nrm{E_j\frac{1}{n^2}\sum_{i\in I}x_i}\leq n\nrm{\frac{1}{n^2}\sum_{i\in I}x_i}\leq 1$$

\nd We can assume that $(E_j)_{j=1}^n$ is a blocking of $( \supp\
x_i)_{i\not\in I}$ and hence
$$\sum_{j=1}^n\nrm{E_j\frac{1}{n^2}\sum_{i\not\in I}x_i}\leq 1$$
which ends the proof of the Lemma.

\

Now we finish the proof of the first inductive step. Applying bounded distortability of
$Y$ and standard diagonalization pick a block subspace $Z$ of $Y$ such that
$\nrm{y}_n/\nrm{y}\leq D\nrm{z}_n/\nrm{z}$ for any non-zero $y,z\in Z_n$, $n\in\N$.

Fix $n\in\N$ and take $x\in Z_n$ as in the Lemma \ref{lem1}. By
the choice of $Z$ for any $y\in Z_n$ we have $\nrm{y}_n\leq
4D\nrm{y}$. By definition of norms $\nrm{\cdot}_n$ it follows that
$$\nrm{E_1y}+\dots+\nrm{E_ny}\leq 4D\nrm{y}$$
for any $y\in Z$ and $n\leq E_1<\dots<E_n$ intervals, $n\in\N$,
which shows that $Z$ is $\ell_1^1$-asymptotic.

\

\nd\textsc{Second inductive step}

Take $1<\al<\omega_1$ and assume now that the theorem holds true
for all $\xi<\al$. If $\al=\beta+1$ for some $\beta<\omega_1$ then
by inductive hypothesis there is a block subspace $W$ of $Y$ which
is $\ell_1^\beta$-asymptotic, and hence $\ell_1^\al$-asymptotic.

If $\al$ is a limit ordinal pick $(\al_n)_n$ with
$\al_n\nearrow\al$ used in the definition of $\cs_\al$. By the
inductive hypothesis we can pick a subspace $W$ of $Y$ such that
$W$ is $\ell_1^{\al_n}$-asymptotic with some constant $C_n\geq 1$
for any $n\in\N$.

Let $\abs{\ \cdot\ }_n$, $n\in\N$, denote the norm given by
$\ell_1^{\al_n}$-asymptoticity of $W$ (Def. \ref{norm}). As before
we will use some special vectors in order to "glue" the original
norm and the new norms. Those vectors - so called special convex
combination, introduced in \cite{ad} - form the crucial tool in
studying properties in particular of mixed and modified mixed
Tsirelson spaces. In order to construct the vectors on
$\ell_1$-$K$-block trees we will slightly generalize the reasoning
from Lemma 4, \cite{klmt} (cf. also Lemma 4.9, \cite{adm}).

\bl\label{lem2} Let $U$ be a Banach space with a basis. Fix
$\eta<\xi<\omega_1$ and assume that $U$ is $\ell_1^1$-asymptotic
with a constant $C_1$, $\ell_1^\eta$-asymptotic and
$I_b(U,K)>\omega^\xi$ for some $K\geq 1$. Then there is a vector
$x\in U$ satisfying $1/K\leq\nrm{x}\leq \abs{x}_\eta\leq 2C_1$.\el

The important part of the Lemma is the fact that the estimates of
the norms in the assertion do not depend on the
$\ell_1^\eta$-asymptoticity constant.

\

\nd\textit{Proof of Lemma \ref{lem2}}. Let $U$ be
$\ell_1^\eta$-asymptotic with constant $C_2$. Let $\ct$ be a
$\ell_1$-$K$-block tree on $U$ with $o(\ct)>\omega^\xi$. We can
assume that for any $(x_i)\in \ct$ also any subsequence
$(x_{i_m})\in \ct$. Put
$$\begin{array}{rl}\cf=& \{(m_1,\dots,m_l)\subset\N:\ m_i\geq \max(supp\ x_i), \ 1\leq i\leq l,
\\ & \ \ \mathrm{for}\ \ \mathrm{some}\ \ (x_1,\dots,x_l)\in\ct\}\end{array}$$
The family $\cf$ is hereditary and either non-compact or, by Prop. 13, \cite{lt1},
compact with $\iota(\cf)\geq o(\ct)>\omega^\xi=\iota(\cs_\xi)$. Hence by Theorem 1.1,
\cite{g1}, there is an infinite $M\subset\N$ with
$$\cs_\xi\cap [M]^{<\infty}\subset\cf$$
 Using Prop. 3.6, \cite{otw}, we get $F\in \cs_\xi\cap [M]^{<\infty}$ and positive scalars
 $(a_m)_{m\in F}$ such that $\sum_{m\in F}a_m=1$ and $\sum_{m\in G}a_m<1/C_2$ for any $G\in
 \cs_\eta$ with $G\subset F$. By definition of $\cf$ there is $(x_i)\in \ct$ such
 that $F=(m_1,\dots,m_l)$ with $m_i\geq \max(supp\ x_i)$ for $1\leq i\leq l$. Let
$x=\sum_{m_i\in F}a_{m_i}x_i$. Since $(x_i)\in\ct$, we have $\nrm{x}\geq 1/K$.

Take now any $\eta$-admissible sequence $E_1<\dots<E_k$. Put
$$J=\left\{j\in\{1,\dots,k\}:\ \min E_j\in \supp\ x_{i_j}\ \mathrm{for\ some}\ i_j\right\}.$$
%$$J=\left\{j\in\{1,\dots,k\}:\ \min E_j\in [\min(supp\ x_{i_j}), \max(supp\ x_{i_j})]\ \mathrm{for\ some}\ i_j\right\}.$$
Let $I=\{ i_j:\ j\in J\}$ and split the sum of the norms as follows
$$\sum_{j=1}^k\nrm{E_jx}\leq \sum_{j=1}^k\nrm{E_j\sum_{i\in I}a_{m_i}x_i}+\sum_{j=1
}^k\nrm{E_j\sum_{i\not\in I}a_{m_i}x_i}$$
 In order to estimate the first part of the sum notice that $G=\{m_{i_j}:\ j\in J\}$ belongs to $\cs_\eta$
since $(\min E_j)_{j=1}^k\in \cs_\eta$ and $m_{i_j}\geq\min E_j$
for any $j\in J$. Hence by $\ell_1^\eta$-asymptoticity of $U$ and
the choice of scalars $(a_m)$ we have
$$\sum_{j=1}^k\nrm{E_j\sum_{i\in I}a_{m_i}x_i}\leq C_2\nrm{\sum_{i\in I}a_{m_i}x_i}\leq C_2\sum_{j\in J}a_{m_{i_j}}\nrm{x_{i_j}}=C_2\sum_{m\in G}a_m\leq 1$$
On the other hand notice that for any $i\not\in I$ and
$j=1,\dots,k$ we have $\min E_j<\min (supp\ x_i)$ whenever $\supp\
x_i\cap E_j\neq\emptyset$. Therefore for sets $J_i=\{j:\ E_j\cap
\supp\ x_i\neq\emptyset\}$, $i\not\in I$, we have $\sharp J_i<\min
(supp\ x_i)$. Hence for any $i\not\in I$ the sequence $(E_j\cap
\supp\ x_i)_{j\in J_i}$ is $1$-admissible and thus
$$\sum_{j=1}^k\nrm{E_j\sum_{i\not\in I}a_{m_i}x_i}\leq\sum_{i\not\in I}a_{m_i}\sum_{j\in J_i}\nrm{E_jx_i}\leq\sum_{i\not\in I}a_{m_i}C_1\nrm{x_i}\leq C_1$$
Putting those estimates together we obtain
$\sum_{j=1}^k\nrm{E_jx}\leq 2C_1$ which ends the proof of Lemma.

\

Now we return to the proof of the second inductive step. Take a
block subspace $Z$ of $W$ such that $\abs{y}_n/\nrm{y}\leq
D\abs{z}_n/\nrm{z}$ for any non-zero $y,z$ in $Z_n$, $n\in\N$.

Since $I_b(Z)>\omega^\al$, then by Lemma 5.8, \cite{jo}, $I_b(Z_n,K)\geq \omega^\al$ for
some $K\geq 1$ and any $n\in\N$.

Fix $n\in\N$. Thus we can use Lemma \ref{lem2} for $U=Z_n$,
$\eta=\al_n$, $\xi=\al_{n+1}$ getting a vector $x\in Z_n$ with
$1/K\leq\nrm{x}$ and $\abs{x}_n\leq 2C_1$. Therefore, by the
choice of $Z$, $\abs{y}_n\leq 2C_1KD\nrm{y}$ for any $y\in Z_n$.
Hence, by the definition of the norms $\abs{\ \cdot\ }_n$,
$$\nrm{E_1y}+\dots+\nrm{E_ky}\leq 2C_1KD\nrm{y}, \ \ \ y\in Z$$
for any $\al_n$-admissible sequence $n\leq E_1<\dots<E_k$ and any
$n\in\N$, i.e. for any $\al$-admissible sequence $E_1<\dots<E_k$,
which shows that $Z$ is $\ell_1^\al$-asymptotic.

\

The second part of the Theorem can be proved in the same way,
replacing $\al$-admissible sequences by $\al$-allowable sequences.
The first inductive step is provided by the assumptions on the
space, the second inductive step follows analogously since Lemma
remains true if one changes "asymptotic" to "strongly asymptotic"
and "admissible sequence" to "allowable sequence". Indeed in the
proof we were working only on $(\min E_j)$, the fact that $E_j$'s
are successive was used only when applying suitable asymptoticity
of the space.

\br As it was mentioned before, the behavior of the norms applied
in the proof was used previously in the study of arbitrary
distortable spaces:

\bnum \item[(i)] Norms $(\nrm{\cdot}_n)_n$ used in the first step
of the proof of Theorem \ref{main} were used in the proof of
arbitrary distortion of the Schlumprecht space \cite{s}. It is
known that such norms give (2-$\vep$)-distortion of Tsirelson
space $T[\cs_1,1/2]$ for any $\vep>0$.

\item[(ii)] In \cite{ot} it was proved that in case of Tsirelson
space $T=T(\cs_1, 1/2)$ the norms $(\nrm{\cdot}_n)_n$ given by
$\ell_1^n$-asymptoticity of $T$ do not arbitrary distort $T$
(\cite{ot}, Thm 2.1, Prop. 1.1). Recall that the block
$\ell_1$-Bourgain index of any infinite dimensional subspace of
$T$ is $\omega^\omega$.

\item[(iii)] In case of the mixed and modified mixed Tsirelson
spaces $T[(\cs_{\al_n}, \theta_n)_n]$ and $T_M[(\cs_{\al_n},
\theta_n)_n]$ studied in \cite{ad,adm,lt2},  the norms
$(\nrm{\cdot}_{\al_n})_n$ distort the whole spaces under certain
conditions on $(\al_n,\theta_n)_n$. In \cite{ad,adm} the special
convex combinations, which we used in our proof, are applied to
produce an asymptotic biorthogonal system.

%\item[(iv)] In case of $T[S_\al,1/2]$ the space is
%(2-$\vep$)-distortable for any $\vep>0$, the proof for the
%original Tsirelson space can be generalized using norms
%$(\nrm{\cdot}_{\al+n})_n$.

\enum \er

\section{The $c_0$ case}

We can formulate in a obvious way analogous definitions of the
$c_0$-block index, denoted here by $J_b(X)$, $c_0^\al$-spreading
models and $c_0^\al$-asymptotic spaces obtaining different
measures of "proximity" of a Banach space to $c_0$. The notion of
$c_0$-block index was investigated in particular in \cite{jo}, in
\cite{g2} higher order $c_0$-spreading models were used in
construction of a strictly singular operator on reflexive
$\ell_1$-asymptotic HI space.

We will sketch here briefly the variant of the reasoning presented in the previous
section, proving the Theorem \ref{main}, and thus in particular Corollary \ref{col}, in
$c_0$ case.

\bt\label{maindual} Let $X$ be a Banach space with a basis. Fix
$\al<\omega_1$. Assume that $J_b(Y)>\omega^\al$ for any block
subspace $Y$ of $X$. Then $X$ contains either an arbitrary
distortable subspace or a $c_0^\al$-asymptotic subspace.

If, additionally, $X$ is $c_0^1$-strongly asymptotic, then $X$
contains either an arbitrary distortable subspace or an
$c_0^\al$-strongly asymptotic subspace.\et

\bp

We shall need suitable norms reflecting $c_0$-asymptoticity of a
space.

\bd\label{normdual} Let $U$ be a Banach space with a basis. Fix
$\al<\omega_1$ and assume $U$ is $c_0^\al$-asymptotic with
constant $C$. The associated norm $\abs{\ \cdot\ }_\al$ is given
by $\abs{x}_\al=\sup\{\abs{\phi(x)}:\ \phi\in U^*,
\abd{\phi}_\al\leq 1\}$, $x\in U$, where
$$\abd{\phi}_\al=\sup\left\{\sum_{i=1}^k\nrd{\phi|_{E_jU}}:\ E_1<\dots<E_k\ \al-\mathrm{admissible}\right\},\ \phi\in
U^*$$ As before, if $U$ is $c_0^\al$-strongly asymptotic, in the
definition of associated norm we use allowable sequences instead
of admissible ones.\ed

\br\label{rem} Clearly $\nrd{\cdot}\leq\abd{\ \cdot\ }_\al\leq
C\nrd{\cdot}$ and $\abs{x}_\al\leq\max\{\nrm{E_jx}:\ 1\leq j\leq
k\}$ for any $x\in U$ and $\al$-admissible $E_1<\dots<E_k$. \er

As before, we will restrict the consideration to the case where $X$ has a block subspace
$Y$ with $D$-bounded distortion.

\

\nd\textsc{First inductive step}

Define on $Y^*$ new equivalent norms as follows:
$$\nrd{\phi}_n=\sup\left\{\sum_{j=1}^n \nrd{\phi|_{E_jY}}:\ E_1<\dots<E_n\ \mathrm{intervals}\right\},\ \
\phi\in Y^*,\ \ n\in\N$$
 Let $\nrm{y}_n=\sup \{\abs{\phi(y)}:\ \phi\in Y^*,\
\nrm{\phi}^*_n\leq 1\}$, $y\in Y$, $n\in\N$.

\bl\label{lem3} Let $U$ be a Banach space with a bimonotone basis.
Fix $n\in\N$ and assume $J_b(U)>\omega$. Then there is a vector
$x\in U$ with $1/2\leq\nrm{x}_n\leq\nrm{x}\leq 2$.\el

\nd \textit{Proof of Lemma \ref{lem3}}. Take a normalized block
sequence $(x_i)_{i=1}^{n^2}\subset U$ which is 2-equivalent to the
unit vector basis of $c_0$ space of dimension $n^2$ and put
$x=\sum_{i=1}^{n^2}x_i$. Obviously $\nrm{x}\leq 2$. Since the
basis is bimonotone we can take normalized functionals
$(\phi_i)_{i=1}^{n^2}\subset U^*$ with $\phi_i(x_i)=1$,
$\phi_i(y)=0$ for any $y\in U$ with $supp\ y\cap \supp\
x_i=\emptyset$, $1\leq i\leq n^2$. Put
$\phi=\frac{1}{n^2}\sum_{i=1}^{n^2}\phi_i$. Since $\phi(x)=1$ it
is enough to show that $\nrd{\phi}_n\leq 2$. Take any
$E_1<\dots<E_n$, define the set $I$ exactly as in the proof of
Lemma \ref{lem1} and proceed computing norms of $\phi_i|_{E_jU}$
instead of
 $E_jx_i$. The second estimate follows from the fact that by the choice of $(\phi_i)$ for any $i\not\in I$ there
is at most one $1\leq j\leq n$ with $\phi_i|_{E_jU}\not\equiv 0$.

\

Coming back to the proof of the first inductive step take a block
subspace $Z$ of $Y$ such that $\nrm{y}_n/\nrm{y}\leq
D\nrm{z}_n/\nrm{z}$ for any non-zero $y,z\in Z_n$ and any
$n\in\N$. For a fixed $n\in\N$ take $x\in Z_n$ as in the Lemma
\ref{lem3} and get $4D\nrm{y}_n\geq \nrm{y}$ for any $y\in Z_n$.
It follows that
$$\nrm{y}\leq 4D\max_{j=1,\dots,n}\nrm{E_jy}$$
for any $y\in Y$ and $n\leq E_1<\dots<E_n$, $n\in\N$, which shows
that $Z$ is $c_0^1$-asymptotic.

\

\nd\textsc{Second inductive step}

Take $1<\al<\omega_1$ and assume that the theorem holds true for
all $\xi<\al$. If $\al=\beta+1$ for some $\beta<\omega_1$ then by
inductive hypothesis there is a block subspace $W$ of $Y$ which is
$c_0^\beta$-asymptotic, and hence also $c_0^\al$-asymptotic.

If $\al$ is a limit ordinal take $(\al_n)_n$ with
$\al_n\nearrow\al$ used in the definition of $\cs_\al$. By the
inductive hypothesis we can pick a subspace $W$ of $Y$ such that
$W$ is $c_0^{\al_n}$-asymptotic with some constant $C_n\geq 1$ for
any $n\in\N$.

Let $\abs{\ \cdot\ }_n$ and $\abd{\ \cdot\ }_n$ denote the norms
on spaces $W$ and $W^*$ given by the $c_0^{\al_n}$-asymptoticity
of $W$. We need the following analogon of Lemma \ref{lem2}:

\bl\label{lem4} Let $U$ be a Banach space with a bimonotone basis.
Fix ordinals $\eta<\xi<\omega_1$ and assume that $U$ is
$c_0^1$-asymptotic with a constant $C_1$, $c_0^\eta$-asymptotic
and $J_b(U,K)>\omega^\xi$ for some $K\geq 1$. Then there is a
vector $x\in U$ with $1/2C_1\leq \abs{x}_\eta\leq\nrm{x}\leq
K$.\el

\nd\textit{Proof of Lemma \ref{lem4}}. Let $U$ be
$c_0^\eta$-asymptotic with constant $C_2$. Proceed as in the proof
of Lemma \ref{lem2}, obtaining normalized block sequence
$(x_i)_{i=1}^l$ in a $c_0$-$K$-block tree, a set
$F=(m_1,\dots,m_l)\in \cs_\xi$ with $m_i\geq \max(supp\ x_i)$,
$1\leq i\leq l$ and suitable positive scalars $(a_m)_{m\in F}$.

Pick normalized functionals $(\phi_i)_{i=1}^l\subset U^*$ with
$\phi_i(x_i)=1$, $\phi_i(y)=0$ for any $y\in U$ with $supp\ y\cap
\supp\ x_i=\emptyset$, $1\leq i\leq l$. Put $x=\sum_{i=1}^lx_i$
and $\phi=\sum_{i=1}^la_{m_i}\phi_i$. Then $\nrm{x}\leq K$. Since
$\phi(x)=1$ it is enough to show that $\abd{\phi}_\eta\leq 2C_1$.

Take any $\eta$-admissible sequence $E_1<\dots<E_k$, define the
sets $J$, $I$ and $J_i$, for $i\not\in I$, exactly as in the proof
of Lemma \ref{lem2} and proceed computing norms of
$\phi_i|_{E_jU}$ instead of $E_jx_i$. To estimate
$\sum_{j=1}^k\nrd{(\sum_{i\in I}a_{m_i}\phi_i)|_{E_jU})}$ use
Remark \ref{rem}. To estimate $\sum_{j=1}^k\nrd{(\sum_{i\not\in
I}a_{m_i}\phi_i)|_{E_jU})}$ use Remark \ref{rem} and the fact,
that $\phi_i|_{E_jU}\equiv 0$ whenever $j\not\in J_i$.

\

Now we return to the proof of the second inductive step. Take a
block subspace $Z$ of $W$ such that $\abs{y}_n/\nrm{y}\leq
D\abs{z}_n/\nrm{z}$ for any non-zero $y,z\in Z_n$, $n\in\N$.

Since $J_b(Z)>\omega^\al$, then $J_b(Z_n,K)\geq \omega^\al$ for
some $K\geq 1$ and any $n\in\N$ (it is easy to check that Lemma
5.8 \cite{jo} is valid also in $c_0$ case). Fix $n\in\N$ and use
Lemma \ref{lem4} for $Z_n$, $\al_n$, $\al_{n+1}$ getting $x\in
Z_n$ with $1/2C_1\leq\abs{x}_n\leq \nrm{x}\leq K$. It follows that
$\nrm{y}\leq 2C_1KD \abs{y}_n$ for any $y\in Z_n$ and thus
$$\nrm{y}\leq 2C_1KD\max_{j=1,\dots,k}\nrm{E_jy}$$
for any $y\in Y$ and $\al$-admissible $E_1<\dots<E_k$, hence $Z$
is $c_0^\al$-asymptotic.

The part for $c_0^\al$-strongly asymptotic spaces follows easily
as in the $\ell_1$ case.\ep

\

\nd Institute of Mathematics\\
Jagiellonian University\\
Reymonta 4, 30-059 Krakow, Poland \\
E-mail: anna.pelczar@im.uj.edu.pl

\end{document}